\documentclass[a4paper,11pt,oneside]{article}
\usepackage[english]{babel}
\usepackage{amsmath,amsfonts,latexsym,amssymb,amsthm}

\begin{document}
\begin{center}{\Large\bf On the largest prime factor of $x^2-1$}\\

\vspace{1cm}
Florian Luca and Filip Najman
\end{center}

\begin{abstract} In this paper, we find all integers $x$ such that $x^{2}-1$ has only prime factors smaller than $100$. This gives some interesting numerical corollaries. For example, for any positive integer $n$ we can find the largest positive integer $x$ such that all prime factors of each of $x, x+1,\ldots, x+n$ are less than 100.
\end{abstract}


\section{Introduction}
\label{sec:1}

For any integer $n$ we let $P(n)$ be the largest prime factor of $n$ with the convention $P(0)=P(\pm 1)=1$. Our main result in this paper is the determination of all integers $x$
satisfying the inequality
 \begin{equation}
\label{eq:1}
P(x^2-1)<100.
\end{equation} 
We also give some interesting corollaries to our main result.  For simplicity, we will consider only solutions with positive values of $x$. 

\medskip

Before stating the results, some history. In 1964, Lehmer \cite{leh} found all positive integer solutions $x$ to the inequality $P(x(x+1))\le 41$. Notice that 
this amounts to finding all odd positive integers $y=2x+1$ such that $P(y^2-1)\le 41$. There are $869$ such solutions. In \cite{Lu}, the first author
found all positive integer solutions of the similar looking equation  $P(x^2+1)<100$. There are $156$ of them.

\medskip

In both \cite{leh} and \cite{Lu}, the method of attack on this question was the following. Assume that $x$ is a positive integers such that 
$P(x^2\pm 1)\le K$, where $K=41$ if the sign is $-$ and $K=100$ if the sign is $+$. Then we can write
\begin{equation}
\label{eq:2}
x^2\pm 1=dy^2,
\end{equation}
where $d$ is squarefree, and $P(dy)\le K$. When  the sign is $-$, then there are $13$ possible primes $p\le 41$ which can participate in the factorization of $d$.
When the sign is $+$, then $p\le 100$ but either $p=2$ or $p$ is a prime which is congruent to $1$ modulo $4$. There is a total of $12$ such primes. 
In both cases, we can write equation \eqref{eq:2} in the form
$$
x^2-dy^2=\mp 1.
$$
Thus, our possible values for $x$ appear as the first coordinate of one of the solutions of at most $2^{13}-1=8191$ Pell equations. For a given Pell equation, 
the sequence $(y_n/y_1)_{n\ge 1}$  forms a Lucas sequence with real roots. The {\it Primitive Divisor Theorem} for Lucas sequences with real roots (see, for example, \cite{car}, or the more general result from \cite{BHV} which applies to all Lucas sequences) says that if $n>6$, then $y_n$ has a prime factor which 
is at least as large as $n-1$. In particular, for the cases treated in \cite{leh} (and \cite{Lu}), it suffices to check the first $42$ (respectively $100$) values 
of the component $x$ of the Pell equations involved there and among these one finds  all possible solutions of the equations considered there.

\medskip  

We follow the same approach in the present paper. The first step is, as previously, to determine the form of all the possible solutions. This is done as explained above via the Primitive Divisor Theorem for the second component of the solutions to Pell equations. Since there are $25$ primes $p<100$, this leads to first solving $2^{25}-1$ Pell equations $x^2-dy^2=1$, the largest one of them having $d$ with $37$ decimal digits. This is clearly impossible using standard algorithms like the continued fractions, as the computations would be too slow and the fundamental solutions of some of the involved Pell equations would have hundreds of millions of decimal digits. Instead, the Pell equations are solved by first computing the regulator of the ring of integers of the corresponding quadratic field, and then from the regulator obtaining a compact representation of the fundamental solution of the Pell equation. The only algorithm known fast enough to compute the huge amount of regulators needed is Buchmann's subexponential algorithm. The output of this algorithm gives exactly the regulator under the Generalized Riemann Hypothesis and unconditionally is only a multiple of the regulator. Next, we check which of the solutions to the Pell equations will lead to solutions to equation \eqref{eq:1}. Finally, a check is performed that proves our search misses no solutions, thus eliminating the apparent dependence of our results on the Generalized Riemann Hypothesis. 

\medskip

{\bf Acknowledgements.}  Part of this work was done while F.~L. visited the Mathematics Department of
the University of Zagreb in February, 2009. He thanks the people of this
Department for their hospitality. Both authors thank Professor Andrej Dujella
and the anonymous referee for useful suggestions. F.~L. was also supported in part by Grants SEP-CONACyT 79685 and PAPIIT 100508. F.~N. was supported by the Ministry of Science, Education and Sports, Republic of Croatia, Grant 037-0372781-2821.

\section{Application of the Primitive Divisor Theorem}
\label{sec:2}

Here, we explain in more detail the applicability of the Primitive Divisor Theorem alluded to in Section
\ref{sec:1} to our problem. 

Let $x$ be an integer such that $x^2-1$ is a product of only the primes up to $97$. We can then write $x^2-1=dy^2$, or, equivalently, $x^2-dy^2=1$, where 
$$
d=2^{a_1}\cdot 3^{a_2}\cdots 97^{a_{25}},\ a_i\in \{0,1\}\text{ for }i=1, \ldots, 25.
$$
The only restriction for $d$ above is that not all the $a_i$s can be zero. This is a Pell equation, so $(x,y)=(x_n,y_n)$, where $x_n+y_n\sqrt d= (x_1+y_1\sqrt d)^n$ for some positive integer $n$ and $x_1+y_1\sqrt d$ is the fundamental solution. As we have $2^{25}-1$ possibilities for $d$, to get a finite number of possible solutions for our equation, we need an upper bound for $n$. 

Let $d$ be fixed, $\eta=x_1+y_1\sqrt d$ and $\zeta=x_1-y_1\sqrt d$. It can be easily seen that 
$$
\left\{
\begin{matrix}
x_n+y_n\sqrt d & = & (x_1+y_1\sqrt d)^n=\eta^n,\\
x_n-y_n\sqrt d & = & (x_1-y_1\sqrt d)^n=\zeta^n.
\end{matrix}
\right.
$$
From here, we deduce that
$$x_n=\frac{\eta^n+\zeta^n}{2}\qquad \text{ and } \qquad y_n=\frac{\eta^n-\zeta^n}{2\sqrt d}.$$
It is easy to see that $y_1$ divides $y_n$. We define
$u_n=y_n/y_1$. As $u_n=(\eta^n-\zeta^n)/(\eta-\zeta)$, it follows that $(u_n)_{n\geq 0}$ is a Lucas sequence of the first kind with real roots $\eta$ and $\zeta$. By a result of Carmichael (see \cite{car}) known as the Primitive Divisor Theorem, it follows that if $n>12$, then $u_n$ has a \emph{primitive divisor} $p$. This primitive divisor has several particular properties, the most important one for us being that 
it satisfies the congruence $p\equiv\pm 1 \pmod n$. This implies that for $n>98$, there exists a prime $p\geq 101$ dividing $u_n$. Thus, $n\leq 98$. 

\section{The algorithm}
\label{sec:3}

We now give our algorithm. Let
$$
S=\{2,3,5,7,11,13,17,19,23,29,31,37,\ldots,97\},
$$ 
be the set of all primes $p\le 100$. The set $S$ has $25$ elements. Let
$R_d$ be the regulator of the ring of integers of the quadratic field $\mathbb Q(\sqrt d)$ and let $x_1+y_1 \sqrt d$ be the fundamental solution of the Pell equation $x^2-dy^2=1$.

\medskip

\par\noindent {\sf
For all $D\in \mathcal P(S)$\\
$\{$\\
$1.\ d=\prod_{p\in D} p$\\
$2.$ Compute $mR_d$ \\
$3.$ Compute a compact representation of $x_m+y_m\sqrt d$\\
$4.$ For $i=1$ to $25$ compute ${\rm ord}_{p_i}y_m$\\
$5.\ z=2^{{\rm ord}_{p_1}y_m}\cdot \ldots 97^{{\rm ord}_{p_{25}}y_m}$\\
$6.$ Compute all the convergents ${p_n}/{q_n}$ of the continued fraction expansion of $\sqrt d$ having $q_n<z$, and check whether $p_n^2-dq_n^2=1$\\
$7.$ If $R_d-\log 2-\log \sqrt d \approx \log z$\\
$\{$\\
$m=1$ and $x_1$ is a solution\\
$8.$ For $i=2$ to $98$\\
$\{$\\
For $j=1$ to $25$ compute ${\rm ord}_{p_j}y_i$\\
$z=2^{{\rm ord}_{p_1}y_i}\cdot \ldots 97^{{\rm ord}_{p_{25}}y_i}$\\
If $i\cdot R_d-\log 2-\log \sqrt d\approx \log z$, $x_i$ is a solution\\
$\}$\\
$\}$\\
$\}$\\
}\\
The algorithm searches through all $2^{25}-1$ possible $d$. 

In step $1$, a value for $d$ is chosen.

In step 2, Buchmann's subexponential algorithm is used to compute $R_d$.
This algorithm returns a multiple $mR_d$ of the regulator $R_d$, unconditionally . If the Generalized Riemann Hypothesis is true, then $m=1$. We will remove the dependence of our algorithm on the Generalized Riemann Hypothesis in step 6.

A \emph{compact representation} of an algebraic number $\beta \in \mathbb Q (\sqrt d)$ is a representation of $\beta$ of the form
\begin{equation}
\beta=\prod_{j=1}^k \left(\frac{\alpha_j}{d_j}\right)^{2^{k-j}}, 
\label{cr}
\end{equation}
where 
$d_j\in
\mathbb Z,\ \alpha_j={(a_j+b_j\sqrt d)}/2 \in \mathbb Q(\sqrt
d),\ a_j,b_j\in \mathbb Z,\ j=1,\ldots, k$, and $k$, $\alpha$ and $d_j$ have $O(\log \sqrt d)$ digits.

In step 3, a compact representation of $x_m+y_m\sqrt d$ is constructed using the algorithms described in \cite{mm}. The reason for using compact representations is that the standard representation of the fundamental solution of the Pell equation has $O(\sqrt d)$ digits. Using the standard representation would make these computations impossible. More details about compact representations can be found in \cite{jw}. 

Since we only have the compact representation of the $x_m+y_m\sqrt d$, in step 4, we use algorithms from \cite{fn} to perform modular arithmetic on the compact representation. The $p$-adic valuations of $y_m$ when $p$ is one of the first 25 primes are computed in the following way. We define $v=\prod_{i=1}^{25}p_i$. We first compute $y_m \pmod v$. If $y_m\equiv 0 \pmod v$, we conclude that $y_m$ is divisible by all 25 primes. Next, for all primes that we now know divide $y_m$, we compute $y_m\pmod{p_i^{c}}$, where $c$ is a sufficiently large constant. We first take $c=15$ and if $y_m\equiv 0\pmod{p_i^{c}}$ still, we then replace $c$ by $2c$ and repeat the computation.

In step 5, $z$, which is defined to be the part of $y_m$ divisible by the first 25 primes, is computed. 

The purpose of step 6 is to remove the dependence of this algorithm on the Generalized Riemann Hypothesis. If the Generalized Riemann Hypothesis is false, it is possible that, without this check, we could miss some solutions in our search. Suppose therefore that $m>1$ and that $y_k$ is such a solution that we are missing, meaning that $y_k$ is divisible only by some of the first 25 primes for some $k$ that is not a multiple of $m$.  But $y_1\mid y_k$, so $y_1$ is divisible only by some of the first 25 primes. Also $y_1\mid y_m$, meaning that $y_1$ divides the part of $y_m$ that is divisible by the first 25 primes, which in our notation is $z$. In other words, $y_1 \mid z$, implying $y_1<z$. As $y_1$ has to be the denominator of a convergent of the continued fraction expansion of $\sqrt d$, it follows that by checking that  the relation $p_n^2-dq_n^2\neq 1$ holds for all $n$ satisfying $q_n<z$, we arrive at a contradiction. This proves that either $m=1$, or that $y_k$ has a prime factor larger than 100 for all positive integers $k$. This implies that our algorithm finds all solutions to equation \eqref{eq:1} unconditionally, and if $y_m=z$, then $m=1$. In practice, $z$ will be a relatively small number, so usually only about $10$ convergents will need to be computed.

Since $x_1\approx y_1 \sqrt d$, it follows that 
$$
R_d=\log(x_1+y_1\sqrt d)\approx \log(2y_1\sqrt d)=\log 2+\log \sqrt d +\log y_1,
$$ 
so we can determine whether $y_1=z$ by the test in step 7. In this test, we took that $a\approx b$ if $|a-b|<0.5$. If $y_1\neq z$, then $y_1\geq 101\cdot z$, so $\log y_1>4.61 +\log z$. This shows that great numerical precision is not needed here. Just in case, we used 10 digits of precision.

Step 8 checks whether any of the $x_n$ are solutions for $n=2$ to $98$. With the purpose of speeding up the algorithm, the fundamental solution was not powered when computing ${\rm ord}_{p_j}y_i$. Instead, we computed the fundamental solution modulo the appropriate integer and then powering modulo that integer. Also, usually not all $n$ have to be checked. This is a consequence of the fact that if a prime $p$ divides $y_k$, then it divides $y_{lk}$ for every positive integer $l$. For example, if we get that $y_2$ is divisible by a prime larger than 100, then $y_4$, $y_6, \ldots$ do not need to be checked.

The running time of the algorithm is dominated by the computation of the regulator in step 2. Step 2 makes one the algorithm run in subexponential time. It is the only part of the algorithm that is not polynomial. The computations were performed on a Intel Xeon E5430. The computation of the regulators took around 12 days of CPU time, while the rest of the computations took around 3 days of CPU time.

Suppose that one wants to find the solutions of $P(x^2-1)<K$ using this algorithm. By the Prime Number Theorem, there are approximately $K/ \log K$ primes up to $K$. This means that the algorithm will loop $2^{K \log K}$ times. The product of all these primes will be of size $O(e^K)$. This means that step 2 will run in $O(e^{\sqrt{K\log K}})$. Thus, the time complexity of the algorithm is $O({\rm exp}(K/\log K +\sqrt{K\log K}))$.

\section{Results}
\label{sec:4}

\newtheorem{tm}{Theorem}
\begin{tm} 

\begin{itemize}

\item[a)] The largest three solutions $x$ of the equation 
$P(x^2-1)<100$
are
$$
x=\left\{
\begin{matrix}
19182937474703818751,\\
332110803172167361,\\
99913980938200001.
\end{matrix}
\right.
$$
\item[b)] The largest solution $x$ of $P(x^4-1)<100$ is $x=4217$.
\item[c)] The largest solution $x$ of $P(x^6-1)<100$ is $x=68$.
\item[d)] The largest $n$ such that $P(x^{2n}-1)<100$ has an integer solution $x>1$ is $n=10$, the solution being $x=2$.
\item[d)] The largest $n$ such that $P(x^{2n}-1)<100$ has an integer solution $x>2$ is $n=6$, the solution being $x=6$.
\item[e)] The equation $P(x^2-1)<100$ has $16167$ solutions.
\item[f)] The greatest power $n$ of the fundamental solution of the Pell equation $(x_1+y_1\sqrt d)^n$ which leads to a solution of our problem is $(2+\sqrt 3)^{18}$; i.e., $n=18$ for $d=3$. The case $d=3$ also gives the most solutions, namely $10$ of them.
\end{itemize}
\end{tm}

\noindent \emph{Proof:}\\
The proof is achieved via a computer search using the algorithm from Section \ref{sec:3}. Part b) is proved by finding the largest square of all the $x$, c) by finding the largest a cube, etc. \qed \\

The largest solution $x$ has $20$ decimal digits, the second largest has $18$ digits, followed by $5$ solutions with $17$ digits and $10$ solutions with $16$ digits. All of the mentioned large solutions are odd. The largest even solution  $x$ has $15$ digits.

\newtheorem{tm2}[tm]{Theorem}
\begin{tm2}
Write the equation \eqref{eq:1} as 
$$
x^2-1=2^{a_1}\cdots97^{a_{25}}.
$$
Then the following hold:
\begin{itemize}
\item[a)] The solution with the largest number of $a_i\neq 0$ is $x=9747977591754401$. For this solution, $17$ of the $a_i$s are non-zero.
\item[b)] The solution with the largest $\sum_{i=1}^{25}{a_i}$ is $x= 19182937474703818751$.  For this solution, $\sum_{i=1}^{25}{a_i}=47$.
\item[c)] The single largest $a_i$ appearing among all solutions is $a_1=27$ and corresponds to the solution $x=4198129205249$.
\end{itemize}
\end{tm2}

\noindent \emph{Proof:}\\
Again, this is done via a computer search using the algorithm from Section \ref{sec:3}. \qed\\

Dabrowski \cite{dab} considered a similar problem, where the prime factors of $x^2-1$ consist of the first $k$ primes $p_1,\ldots, p_k$. He formulated the following conjecture.

\newtheorem{thm54}[tm]{Conjecture}
\begin{thm54}
The Diophantine equation 
$$x^2-1=p_1^{\alpha_1}\cdots p_k^{\alpha_k}$$
has exactly  $28$ solutions $(x;\alpha_1,\ldots,\alpha_k)$ in positive integers, as follows:
\begin{itemize}
\item[a)]$(3;3)$,
\item[b)]$(5;3,1),(7;4,1),(17;5,2)$,
\item[c)]$(11; 3, 1, 1), (19; 3, 2, 1), (31; 6, 1, 1), (49; 5, 1, 2), (161; 6, 4, 1),$
\item[d)]$(29; 3, 1, 1, 1), (41; 4, 1, 1, 1), (71; 4, 2, 1, 1), (251; 3, 2, 3, 1), (449; 7, 2, 2, 1),$
$(4801; 7, 1, 2, 4),(8749; 3, 7, 4, 1),$
\item[e)]$(769; 9, 1, 1, 1, 1), (881; 5, 2, 1, 2, 1), (1079; 4, 3, 1, 2, 1), (6049; 6, 3, 2, 1, 2),$ 
$(19601; 5, 4, 2, 2, 2),$
\item[f)]$(3431; 4, 1, 1, 3, 1, 1), (4159; 7, 3, 1, 1, 1, 1), (246401; 8, 6, 2, 1, 1, 2),$
\item[g)]$(1429; 3, 1, 1, 1, 1, 1, 1), (24751; 5, 2, 3, 1, 1, 1, 1), (388961; 6, 4, 1, 4, 1, 1, 1),$
\item[h)]$(1267111; 4, 3, 1, 1, 3, 1, 1, 2).$
\end{itemize}
\label{dabslut}
\end{thm54}
The main result of \cite{dab} is that Conjecture \ref{dabslut} is true for $k\leq 5$.
From our data, we confirm Dabrowski's conjecture in a much wider range.
\newtheorem{thm56}[tm]{Theorem}
\begin{thm56}
Conjecture \ref{dabslut} is true for $k\leq 25$.
\end{thm56}

\noindent\emph{Proof:}\\
This is done by simply factoring all $x^2-1$, where $x$ is a solution of \eqref{eq:1}.\qed\\

The next theorem follows also trivially from our results. Recall that a positive integer $n$ is {\it $K$-smooth} if $P(n)\le K$. In particular, the main result of our paper is the determination of all the $100$-smooth positive integers of the form $x^2-1$.

\newtheorem{tm3}[tm]{Corollary}
\begin{tm3}
Let $t$ be the largest odd solution, $t=19182937474703818751$ and $s$ be the largest even solution, $s=473599589105798$ of equation \eqref{eq:1}.
\begin{itemize}
\item[a)] The largest consecutive $100$-smooth integers are $x$ and $x+1$ where $x=(t-1)/2$ is the largest solution of $P(x(x+1))<100.$
\item[b)] The largest consecutive even $100$-smooth integers are $t-1$ and $t+1$.
\item[c)] The largest consecutive odd $100$-smooth integers are $s-1$ and $s+1$.
\item[d)] The largest triangular $100$-smooth integer is $(t^2-1)/8$.
\end{itemize}
\end{tm3}

As we mentioned in the Introduction, the problem of finding two consecutive $K$-smooth integers was examined by Lehmer in \cite{leh} in the sixties. At that time, he was able to solve the above problem for  the values $K\le 41$. The advance of both computing power and theoretical arguments (namely, the compact represenations) allow us to solve the much harder problem of finding consecutive $K$-smooth integers for any $K\leq 100$. Note that, as was already mentioned at the end of Section \ref{sec:3}, the difficulty of this problem grows exponentially with $K$.

Our results can also be applied to finding $k$ consecutive $K$-smooth integers, for any integer $k$. 
We obtain the following results. 
\newtheorem{thm57}[tm]{Corollary}
\begin{thm57}
The largest integer $x$ satisfying
$$P(x(x+1)\ldots(x+n))<100,$$
for a given $n$, are given in the following table: 
\begin{center}
\begin{tabular}{|c|c|}
\hline
$n$ & $x$\\
\hline
$1$ & $9591468737351909375$\\
\hline
$2$ & $407498958$\\
\hline
$3$ & $97524$\\
\hline
$4$ & $7565$\\
\hline
$5$ & $7564$\\
\hline
$6$ & $4896$\\
\hline
$7$ & $4895$\\
\hline
$8$ & $284$\\
\hline
\end{tabular}
\end{center}
\end{thm57}

\noindent \emph{Proof:}\\
To find $k$ consecutive $100$ smooth integers, we first create a list of all pairs of consecutive $100$-smooth integers. Every odd solution of our starting problem will give us one such pair. This is because if $x$ is an odd solution, then $(x-1)/2$ and $(x+1)/2$ are consecutive $100$-smooth integers. Once this list is created, we search for overlaps in these pairs and obtain our results. \qed\\

In the recent paper \cite{ST}, Shorey and Tijdeman proved several extensions of some irreducibility theorems due to Schur. The main results of \cite{ST} rely heavily on Lehmer's results from \cite{leh}. Thus, replacing, for example Lemma 2.1 in \cite{ST} by our results, it is likely that the main results from 
\cite{ST} can be extended in a wider range.

\medskip

\textbf{Remark.}
The tables produced by our computations can be found on the webpage {\sf http://web.math.hr/$\sim$fnajman}.

\noindent \small{FLORIAN LUCA}\\
\small{INSTITUTO DE MATEMATICAS,\\ UNIVERSIDAD NACIONAL AUTONOMA DE MEXICO,\\ C.P. 58089, MORELIA, MICHOACAN, MEXICO}\\
\emph{E-mail address:} fluca@matmor.unam.mx
\vspace{0.5cm}\\
\small{FILIP NAJMAN}\\
\small{DEPARTMENT OF MATHEMATICS,\\ UNIVERSITY OF ZAGREB,\\ BIJENI\v CKA CESTA 30, 10000 ZAGREB, CROATIA}\\
\emph{E-mail address:} fnajman@math.hr

\end{document}